\documentclass[12pt]{amsart}

\usepackage[margin=1.5in]{geometry}
\usepackage{amsmath}
\usepackage{amsfonts}
\usepackage{amssymb}

\usepackage{hyperref}

 \newtheorem{lemma}{Lemma}[section]
  \newtheorem{theorem}[lemma]{Theorem}

  \newtheorem{remark}[lemma]{Remark}

\numberwithin{equation}{section}

\begin {document}
\title[Kihara's curve]{The rank and generators of Kihara's elliptic curve with torsion $\mathbb{Z}/4\mathbb{Z}$
 over $\mathbb{Q}(t)$}

\author[A. Dujella,  I. Gusi\'c and P. Tadi\'c]{Andrej Dujella, Ivica Gusi\'c and Petra Tadi\'c}

\address[Andrej Dujella]{Department of Mathematics\\ University of Zagreb\\
Bijeni\v{c}ka cesta 30, 10000 Zagreb\\ Croatia}
\email[Andrej Dujella]{duje@math.hr}

\address[Ivica Gusi\'c]{Faculty of Chemical Engin. and Techn. \\
University of Zagreb \\ Maruli\'cev trg 19, 10000 Zagreb \\ Croatia}
\email[Ivica Gusi\'c]{igusic@fkit.hr}

\address[Petra Tadi\'c]{Marti\' ceva 23, 10000 Zagreb\\ Croatia}
\email[Petra Tadi\'c]{petra.tadic.zg@gmail.com}

\thanks{The authors were supported by the Croatian Science Foundation under the project no. 6422.}

\date{}
\maketitle

\begin{abstract}
For the elliptic curve $E$ over $\mathbb{Q}(t)$ found by Kihara,
with torsion group $\mathbb{Z}/4\mathbb{Z}$
and rank $\geq 5$, which is the current record for the rank of such curves,
by using a suitable injective specialization, we determine exactly the rank and generators of
$E(\mathbb{Q}(t))$.
\end{abstract}

\footnotetext{ {\it 2010 Mathematics Subject Classification.}
11G05, 14H52.\\
{\it Key words and phrases.} Elliptic curve, specialization
homomorphism, rank, torsion, generators, mwrank}

\section{Introduction}  \label{section1}
Let $E$ be an elliptic curve over $ \mathbb{Q}$.
By the Mordell-Weil theorem, the group $E( \mathbb{Q})$ of rational points
on $E$ is a finitely generated abelian group. Hence, it is the
product of the torsion group and $r\geq 0$ copies of the infinite cyclic
group: $E( \mathbb{Q}) \cong E( \mathbb{Q})_{\rm
tors} \times { \mathbb{Z} }^r$.
An important question in arithmetic geometry is to characterize
which groups are possible as $E( \mathbb{Q})$ for an elliptic curve $E$ over $\mathbb{Q}$.
Thus, we may ask which torsion groups are possible, which ranks are possible,
and, most ambitiously, which combinations of torsion and rank are possible.
By Mazur's theorem, we know that $E( \mathbb{Q})_{\rm tors}$ is one
of the following 15 groups:
$ \mathbb{Z}/n\mathbb{Z}$ with $1\leq n \leq 10$ or $n=12$, $ \mathbb{Z}/2\mathbb{Z} \times \mathbb{Z}/2m\mathbb{Z}$ with $1\leq
m\leq 4$.
On the other hand, it is not known which values of rank $r$ are
possible for elliptic curves over $ \mathbb{Q}$. The folklore
conjecture is that the rank can be arbitrary large, but it seems to be
very hard to find examples with high rank. The current record is an
example of an elliptic curve over $\mathbb{Q}$ with rank $\geq 28$,
found by Elkies in 2006.
There is even a stronger conjecture that for any of 15 possible
torsion groups $T$ we have $B(T)=\infty$, where
$B(T)=\sup \{ {\rm rank}\,(E(\mathbb{Q})) \,:\,
\mbox{torsion group of $E$ over $ \mathbb{Q}$ is $T$} \}$.
The current records for ranks for each of the $15$ possible torsion groups
can be found on the web page \cite{D-tors}.

Similar questions can be asked for elliptic curves over the field of rational
functions $\mathbb{Q}(t)$. Again, we have $15$ possible torsion groups,
and we define $G(T)=\sup \{ {\rm
rank}\, E(\mathbb{Q}(t)) \,:\, E( \mathbb{Q}(t))_{\rm tors} \cong T \}$.
In the next table, taken from \cite{D-fam}, we give the current records for
the rank of elliptic curves over $\mathbb{Q}(t)$ with prescribed torsion group.

\medskip

{\footnotesize

\begin{flushleft}
\begin{tabular}{|@{\quad}c@{\quad}|@{\quad}c@{\quad}|@{\quad}l@{\quad}|}
\hline\rule{0pt}{12pt} $T$ & $G(T) \geq $ &
Author(s)\\[2pt]
\hline\rule{0pt}{12pt}
0 &          18 &        Elkies (2006) \\[5pt]
  $\mathbb{Z}/2\mathbb{Z}$  &       11 &       Elkies  (2009) \\[5pt]
  $\mathbb{Z}/3\mathbb{Z}$  &        7 &       Elkies  (2007) \\[5pt]
  $\mathbb{Z}/4\mathbb{Z}$ &         5 &       Kihara (2004),  Elkies (2007)  \\[5pt]
  $\mathbb{Z}/5\mathbb{Z}$  &        3 &       Lecacheux  (2001), Eroshkin (2009), MacLeod (2014) \\[5pt]
  $\mathbb{Z}/6\mathbb{Z}$  &        3 &       Lecacheux  (2001), Kihara (2006), Eroshkin (2008),\\
  & &  Woo (2008), Dujella \& Peral (2012), MacLeod (2014) \\[5pt]
  $\mathbb{Z}/7\mathbb{Z}$  &        1 &       Kulesz (1998), Lecacheux (2003), Rabarison (2008),\\
  & &  Harrache (2009), MacLeod (2014) \\[5pt]
  $\mathbb{Z}/8\mathbb{Z}$  &        2 &       Dujella \& Peral (2012), MacLeod (2013) \\[5pt]
  $\mathbb{Z}/9\mathbb{Z}$  &        0 &       Kubert (1976) \\[5pt]
  $\mathbb{Z}/10\mathbb{Z}$ &        0 &       Kubert (1976) \\[5pt]
  $\mathbb{Z}/12\mathbb{Z}$ &        0 &       Kubert (1976) \\[5pt]
$\mathbb{Z}/2\mathbb{Z} \times \mathbb{Z}/2\mathbb{Z}$ &    7 &       Elkies  (2007) \\[5pt]
$\mathbb{Z}/2\mathbb{Z} \times \mathbb{Z}/4\mathbb{Z}$ &    4 &  Dujella \& Peral (2012) \\[5pt]
$\mathbb{Z}/2\mathbb{Z} \times \mathbb{Z}/6\mathbb{Z}$ &     2 & Dujella \& Peral (2012), MacLeod (2013) \\[5pt]
$\mathbb{Z}/2\mathbb{Z} \times \mathbb{Z}/8\mathbb{Z}$ &     0 &
Kubert (1976)
 \\[2pt]
\hline
\end{tabular}
\end{flushleft}
}%

\medskip

Note that in the previous table only the lower bounds for the rank are given.
Indeed, it seems that only for the torsion group $\mathbb{Z}/2\mathbb{Z} \times \mathbb{Z}/4\mathbb{Z}$
the exact rank over $\mathbb{Q}(t)$ of the record curve can be found in literature.
In fact, in \cite{D-P},
Dujella and Peral proved that the corresponding curve, obtained from the so called Diophantine triples,
has rank equal to $4$ and they provide the generators for the group. The proof uses the method
introduced by Gusi\'c and Tadi\'c in \cite{G-T1} for an efficient search for
injective specializations. In this paper, we will prove an analogous result for the
curve with record rank over $\mathbb{Q}(t)$ with torsion group  $\mathbb{Z}/4\mathbb{Z}$
found by Kihara \cite{Kih2} (see Theorem \ref{th:main}).
Here we will use results from the recent paper \cite{G-T2}, where the authors generalize and
extend their method from \cite{G-T1}. In particular, by results of \cite{G-T2}, now the method
can be applied to curves with only one rational $2$-torsion point.

Our main tool is \cite[Theorem 1.3]{G-T2}. It deals with elliptic curves $E$
given by $y^2 = x^3 + A(t)x^2 + B(t)x$, where $A,B \in \mathbb{Z}[t]$, with exactly
one nontrivial $2$-torsion point over $\mathbb{Q}(t)$.
If $t_0 \in \mathbb{Q}$ satisfies the condition that for every nonconstant
square-free divisor $h$ of $B(t)$ or $A(t)^2 - 4B(t)$ in $\mathbb{Z}[t]$
the rational number $h(t_0)$ is not a square in $\mathbb{Q}$,
then the specialized curve $E_{t_0}$ is elliptic and the specialization homomorphism at
$t_0$ is injective.
If additionally  there exist $P_1,\ldots,P_r\in E(\mathbb{Q}(t))$ such that
$P_1(t_0),\ldots,P_r(t_0)$ are the free generators of $E(t_0)(\mathbb{Q})$, then
$E(\mathbb{Q}(t))$ and $E(t_0)(\mathbb{Q})$ have the same rank $r$, and $P_1,\ldots,P_r$ are the free generators of $E(\mathbb{Q}(t))$.

We can mention here that by the methods from \cite{G-T2}, it is easy to show that
the general families of curves with torsion $\mathbb{Z}/10\mathbb{Z}$, $\mathbb{Z}/12\mathbb{Z}$
and $\mathbb{Z}/2\mathbb{Z} \times \mathbb{Z}/8\mathbb{Z}$ given in Kubert's paper \cite{Kub}
all have rank $0$ over $\mathbb{Q}(t)$, as expected from the corresponding entries in the
above table (see Remark \ref{rem:kubert}).

\section{Kihara's curve with rank $\geq 5$}
In 2004, Kihara \cite{Kih2} constructed a curve over $\mathbb Q(t)$ with
torsion group $\mathbb{Z}/4\mathbb{Z}$ and rank $\geq 5$. This improved his previous
result \cite{Kih} with rank $\geq 4$. We briefly describe Kihara's construction.
The quartic curve $H$ given by the equation
\begin{equation} \label{eq:H}
(x^2-y^2)^2 + 2a(x^2+y^2)+b=0
\end{equation}
is considered. Forcing five points with coordinates of the form $(r,s)$, $(r,u)$, $(s,p)$,
$(u,q)$, $(p,m)$ to satisfy (\ref{eq:H}) leads to a system of certain quadratic
Diophantine equations, for which a parametric solution is found.
By the transformation $X=(a^2-b)y^2/x^2$ and $Y=(a^2-b)y(b+ax^2+ay^2)/x^3$, we get
from $H$ the elliptic curve $E$ with equation
\begin{equation} \label{eq:E}
Y^2 = X(X^2 + (2a^2 + 2b)X+ (a^2 - b)^2).
\end{equation}
We can write (\ref{eq:E}) in the form
$$ Y^2 = X^3 + A(t)X^2 + B(t)X, $$
where
{\small
\begin{eqnarray*}
A(t) &\!\!=\!\!& 26244t^{52}-1454544t^{51}+39753720t^{50}-693511680t^{49}+8373526532t^{48}-71499362784t^{47}\\ & & \hspace*{-10mm}\mbox{}+419498533496t^{46}-1416310911728t^{45}-901962029472t^{44}+45707741530384t^{43}\\
& & \hspace*{-10mm}\mbox{}-311918617874472t^{42}+1112855854012288t^{41}-445015546446924t^{40}-24074707538255904t^{39} \\
& & \hspace*{-10mm}\mbox{}+201749209147935960t^{38}-1090449631523616720t^{37}+4657677879944558620t^{36}\\
& & \hspace*{-10mm}\mbox{}-16789722010744550048t^{35}+52691499185487218416t^{34}-146529244725968515840t^{33}\\
& & \hspace*{-10mm}\mbox{}+365082905946272744264t^{32}-820936801892256346048t^{31}+1674263798355290924784t^{30}\\
& & \hspace*{-10mm}\mbox{}-3107412066279198054752t^{29}+5260774429646277796992t^{28}-8138423149394923126752t^{27}\\
& & \hspace*{-10mm}\mbox{}+11524246769870127828848t^{26}-14973162439865722690560t^{25}+17926531556504852382024t^{24}\\
& & \hspace*{-10mm}\mbox{}-19931464457985164353856t^{23}+20851433797507505280624t^{22}-20915937991381166319520t^{21}\\
& & \hspace*{-10mm}\mbox{}+20534733115474642866140t^{20}-19974922175127737352208t^{19}+19137997688945377556248t^{18}\\
& & \hspace*{-10mm}\mbox{}-17635555889425637775104t^{17}+15136697914178618626228t^{16}-11721796952272669744224t^{15}\\
& & \hspace*{-10mm}\mbox{}+7949459617497409864344t^{14}-4573957484971788140208t^{13}+2138285815612380734304t^{12}\\
& & \hspace*{-10mm}\mbox{}-748888327112027254832t^{11}+152762005274598499512t^{10}+14561852207708951680t^9\\
& & \hspace*{-10mm}\mbox{}-26448842122489903356t^8+11493414948906352480t^7-2992433840971722568t^6\\
& & \hspace*{-10mm}\mbox{}+511890043268967024t^5-57486161177498812t^4+4053242286328704t^3-161507610671616t^2\\
& & \hspace*{-10mm}\mbox{}+2684575014912t+8707129344,
\end{eqnarray*}
\begin{eqnarray*}
 B(t) &\!\!=\!\!& 256t^2(t-1)^4(t+1)^4(t-3)^2(t-5)^2(3t-1)^2(t-2)^4(t^2+2)^2(t^2-2t+3)^2 \\
 & & \hspace*{-10mm}\mbox{}\times (t^2+6t-1)^2(7t^2-18t+23)^2(3t^2-2t+7)^2(t^3+4t^2-5t+16)^2 \\
 & & \hspace*{-10mm}\mbox{}\times (3t^4-17t^3+27t^2-43t+6)^2 (2t^4-7t^3+9t^2-11t-5)^2(2t^4-17t^3+27t^2-25t+1)^2 \\
 & & \hspace*{-10mm}\mbox{}\times (t^4-4t^3+6t^2-12t+1)^2(5t^4-17t^3+27t^2-79t+16)^2(t^4-28t^3+54t^2-92t+41)^2 \\
 & & \hspace*{-10mm}\mbox{}\times (t^4-9t^3+15t^2-19t+4)^2.
\end{eqnarray*}
}%

Kihara in \cite{Kih2} found five independent points $P_1,\ldots, P_5$ on this curve,
corresponding to the five points on $H$ mentioned above,
showing that the rank of $E$ over $\mathbb Q(t)$ is $\geq 5$. The torsion subgroup is $\mathbb Z/4\mathbb Z$. Indeed, the point $T_1=(a^2-b,2a(a^2-b))$ on (\ref{eq:E}) is of order $4$ since $2T_1=(0,0)$ and
$4T_1=\mathcal{O}$.

Our goal is to prove that the rank of $E$ over $\mathbb Q(t)$ is exactly equal to $5$
and to find the generators of $E(\mathbb Q(t))$.
Computations with several specializations
indicate that $P_1,P_2,P_3,P_4,P_5$ are not generators of $E(\mathbb Q(t))$.
Indeed, from our results it will follow that they generate a subgroup of index $16$ in $E(\mathbb Q(t))$.

In fact, it holds that $P_1+P_i\in 2E(\mathbb Q(t))$ for $i=2,3,4,5$,
i.e. there exist points $W_2,W_3,W_4,W_5$ of $E(\mathbb Q(t))$
such that $P_1+P_i=2W_i, \ i=2,3,4,5$.
Since the torsion subgroup is $\mathbb Z/4\mathbb Z$,
there are two choices for each $W_i$. We choose one of them. The $x$-coordinates of these points are
{\small
\begin{eqnarray*}
x(W_2) &\!\!=\!\!&
4(t-3)(3t-1)(t^2+2)(t^4-4t^3+6t^2-12t+1)(3t^4-17t^3+27t^2-43t+6) \\
& & \hspace*{-10mm}\mbox{}\times (t^4-9t^3+15t^2-19t+4)(7t^2-18t+23)(t^4-28t^3+54t^2-92t+41)(t-2)^2 \\
& & \hspace*{-10mm}\mbox{}\times (9t^{12}-110t^{11}+576t^{10}-2333t^9+7802t^8-19832t^7+39488t^6-57374t^5 \\
& & \hspace*{-10mm}\mbox{}+61421t^4-42914t^3+16488t^2-701t-216)^2 (t+1)^4, \\
x(W_3) &\!\!=\!\!&
-16(t-1)^3(t+1)^3(t-3)(t^2+2)(t^2+6t-1)(t^3+4t^2-5t+16) \\
& & \hspace*{-10mm}\mbox{}\times
(t^4-9t^3+15t^2-19t+4)(2t^4-17t^3+27t^2-25t+1)(3t^4-17t^3+27t^2-43t+6) \\
& & \hspace*{-10mm}\mbox{}\times (3t^2-2t+7)(t-2)^2t^2(3t-1)^2
(t-5)^2(t^4-28t^3+54t^2-92t+41)^2 \\
& & \hspace*{-10mm}\mbox{}\times(2t^4-7t^3+9t^2-11t-5)^2, \\
x(W_4) &\!\!=\!\!&
16t(t-1)^3(t+1)(3t-1)(t^2+2)(t^2+6t-1)(t^3+4t^2-5t+16) \\
& & \hspace*{-10mm}\mbox{}\times (3t^4-17t^3+27t^2-43t+6)
(t^4-9t^3+15t^2-19t+4)(2t^4-17t^3+27t^2-25t+1)\\
& & \hspace*{-10mm}\mbox{}\times (3t^2-2t+7)(t-2)^2(t-5)^2(7t^2-18t+23)^2
(t^4-4t^3+6t^2-12t+1)^2 \\
& & \hspace*{-10mm}\mbox{}\times (2t^4-7t^3+9t^2-11t-5)^2
(3t^7-29t^6+51t^5-136t^4+175t^3-267t^2+179t-24)^2 / \\
& & \hspace*{-5mm} (8-55t-29t^2+103t^3-120t^4+69t^5-27t^6+3t^7)^2, \\
x(W_5) &\!\!=\!\!&
16(t-1)^3(t-2)^3(t-3)^2(t^2+2)^2(3t^2-2t+7)^2(t^3+4t^2-5t+16)^2\\
& & \hspace*{-10mm}\mbox{}\times (t^4-28t^3+54t^2-92t+41)
(t^4-4t^3+6t^2-12t+1)(3t^4-17t^3+27t^2-43t+6)\\
& & \hspace*{-10mm}\mbox{}\times (t^4-9t^3+15t^2-19t+4)(5t^4-17t^3+27t^2-79t+16)
(2t^4-17t^3+27t^2-25t+1)\\
& & \hspace*{-10mm}\mbox{}\times (t^2+6t-1)(t-5)(3t-1)(t+1)t.
\end{eqnarray*}
}%
We also give the $x$-coordinate of $P_1$: 
{\small
\begin{eqnarray*}
x(P_1) &\!\!=\!\!&
16t(t-2)^2(t-3)(t-5)(3t-1)(t^2+2)(3t^2-2t+7)(7t^2-18t+23)\\ & & \hspace*{-10mm}\mbox{}\times (t^2+6t-1)(t^3+4t^2-5t+16)(2t^4-17t^3+27t^2-25t+1) \\ 
& & \hspace*{-10mm}\mbox{}\times (2t^4-7t^3+9t^2-11t-5)(3t^4-17t^3+27t^2-43t+6)(t^4-4t^3+6t^2-12t+1)\\
& & \hspace*{-10mm}\mbox{}\times (t^2-2t+3)(5t^4-17t^3+27t^2-79t+16)(t^4-9t^3+15t^2-19t+4)\\
& & \hspace*{-10mm}\mbox{}\times (t^4-28t^3+54t^2-92t+41)
(3t^{13}-128t^{12}+1185t^{11}-5018t^{10}+13628t^9-27704t^8\\
& & \hspace*{-10mm}\mbox{} +44162t^7-63956t^6+84827t^5-100976t^4+92061t^3-52802t^2+10662t-552)^2/ \\
& & \hspace*{-5mm}\mbox{} (12t^{11}-219t^{10}+1699t^9-7248t^8+21004t^7-45434t^6+72862t^5-90128t^4\\
& & \hspace*{-10mm}\mbox{} +77496t^3-46283t^2+10095t-768)^2.
\end{eqnarray*}
}%
The points $P_1,W_2,\ldots,W_5$ are a natural guess for the generators, and we will show that
this is indeed true by proving the following theorem in the next section.

\begin{theorem} \label{th:main}
The elliptic curve $E$ over $\mathbb Q(t)$ has rank equal to $5$ with  free generators
the points $P_1,W_2,W_3,W_4,W_5$ and the torsion group is $\mathbb Z/4\mathbb Z$.
\end{theorem}

\section{An injective specialization}

As described in the introduction, we use \cite[Theorem 1.3]{G-T2} to find rational numbers
$t_0$ for which the specialization map at $t_0$ is injective. The condition is
that for each nonconstant square-free divisor $h$ of $B(t)$ or $A(t)^2-4B(t)$ in $\mathbb Z[t]$
the rational number $h(t_0)$ is not a square in $\mathbb{Q}$. The condition is easy to check,
and we can find many rationals $t_0$ satisfying it. However, the coefficients of the curve $E$
are polynomials with large degrees and coefficients. Thus, for the success of our approach,
it is crucial to find suitable specialization $t_0$ of reasonably small height.
Furthermore, we need a specialization for which the rank of $E_{t_0}$ over $\mathbb{Q}$ is
equal to $5$, so it is reasonable to consider only specializations for which the the root number
of $E_{t_0}$ is $-1$ (conjecturally implying that the rank is odd).

We find that the specialization at $t_0=-\frac {11}4$ satisfies
all requirements. It remains to compute the rank and generators of $E_{-11/4}$.
For that purpose, we use the excellent program \cite{mwrank} of Cremona,
which is included in the program package Sage \cite{sage}.
By extending significantly the default precision (we use options {\tt -p 800 -b 11}),
we get the elliptic curve $E_{-11/4}$ over $\mathbb Q$, given by the equation
{\scriptsize
\begin{eqnarray*}
y^2 &\!\!=\!\!& x^3+484371205173916954475505177386303655600428018856419825361x^2 \\
& & \hspace*{-30mm}\mbox{}+39643973962637685622080358722553196979593524013544183541184924577327568643868932220629361508907165935719959040000x,
\end{eqnarray*}
}%
which is of rank 5 with five free generators $G_1,\ldots,G_5$
and the generator of the torsion group $T_0$, given by their $x$-coordinates
{\tiny
\begin{eqnarray*}
x(T_0) &\!\!=\!\!& -199107945503532544541442922607063711556495848199368164800 \\
x(G_1) &\!\!=\!\!& \frac{35128929795293330966584382924686967322464932483259844000000}{49} \\
x(G_2) &\!\!=\!\!&-\frac{13315895444389669790150691801504854333267656768880200356488633890879100326435093056000}{76159758997263590677307690401}\\
x(G_3) &\!\!=\!\!& -\frac{536936685248255028234168808334606137477987482678373004839681994240}{1416844881} \\
x(G_4) &\!\!=\!\!& \frac{477616878094060797794543416082366648044181864476558959811339468000}{9223369} \\
x(G_5) &\!\!=\!\!& \frac{727946070485627419024954469098052196447100359265007972065657045989202204731022849600}{13386309077372650899951050449}
\end{eqnarray*}
}%
The rank 5 and generators are also confirmed
in the most recent version V2.20-10 of Magma \cite{magma}
(by the function {\tt MordellWeilShaInformation} with
option {\tt SetClassGroupBounds("GRH")}). Now denote by $P_1^*, W_2^*, \ldots, W_5^*$ the points obtained from
$P_1,W_2,\ldots,W_5$ after specialization $t\mapsto -\frac {11}4$.
We easily get that
\begin{eqnarray*}
P_1^* &\!\!=\!\!& T_0+ G_5, \\
W_2^* &\!\!=\!\!& T_0- G_3+ G_4, \\
W_3^* &\!\!=\!\!& T_0+ G_1+ G_2- G_3+ G_4- G_5, \\
W_4^* &\!\!=\!\!& T_0- G_2+ G_3- G_4+ G_5, \\
W_5^* &\!\!=\!\!& T_0+ G_4.
\end{eqnarray*}
Here all points are chosen up to sign of $y$-coordinates.
It is easy to check that the matrix of this base change (modulo torsion)
is of determinant $\pm 1$ (the signs of the determinant depends on the choice of the signs of
$y$-coordinates),
so we see that $P_1^*,W_2^*,\ldots,W_5^*$ are free generators of the elliptic curve
$E_{-11/4}$ over $\mathbb{Q}$. From the comments at the end of the introduction,
we see that \cite[Theorem 1.3]{G-T2} now implies
that $E$ has rank $5$ over $\mathbb Q(t)$ and that $P_1,W_2,W_3,W_4,W_5$ are its free generators. Since $E$ has a point of fourth order and the torsion group of $E_{-11/4}(\mathbb Q)$ is $\mathbb Z/4\mathbb Z$, we conclude that the torsion group of $E(\mathbb Q(t))$ is also $\mathbb Z/4\mathbb Z$.
\qed

\begin{remark} \label{rem:kubert}
{\rm 
Here we prove that the elliptic curves with torsion $\mathbb Z/10\mathbb Z$, $\mathbb Z/12\mathbb Z$, $\mathbb Z/2\mathbb Z\times\mathbb Z/8\mathbb Z$ and with rank $\geq 0$ listed in Kubert's paper \cite{Kub} have rank equal to $0$ over $\mathbb{Q}(t)$.

For torsion group $\mathbb Z/10\mathbb Z$, the curve is given by the equation
$$y^2=x(x^2-(2t^2-2t+1)(4t^4-12t^3+6t^2+2t-1)x+16t^5(t-1)^5(t^2-3t+1)).$$ 
The specialization for $t_0=6$ satisfies the condition of \cite[Theorem 1.3]{G-T2}, 
and the specialized elliptic curve has rank $0$, which proves our claim for that torsion group.

For torsion group $\mathbb Z/12\mathbb Z$, the curve is 
$$y^2=x(x^2+(t^8-12t^6-48t^5-162t^4-480t^3-540t^2-624t-183)x+1024(t^2+3)^2(t+1)^6),$$ 
and the specialization which satisfies the condition of \cite[Theorem 1.3]{G-T2} 
and has rank $0$ is $t_0=11$.  

Finally, for torsion group $\mathbb Z/2\mathbb Z\times\mathbb Z/8\mathbb Z$, 
we have the curve 
$$y^2=x(x+(t^2-1)^4)(x+16t^4).$$ 
The specialization for $t_0=6$ satisfies the condition of \cite[Theorem 1.1]{G-T2}, and the specialized elliptic curve has rank $0$, which proves our claim in that case. 
}
\end{remark}

\end{document}